\documentclass[11pt,a4paper]{article}

\usepackage{graphicx}%
\usepackage{multirow}%
\usepackage{amsmath,amssymb,amsfonts}%
\usepackage{amsthm}%
\usepackage{mathrsfs}%
\usepackage[title]{appendix}%
\usepackage{xcolor}%
\usepackage{textcomp}%
\usepackage{manyfoot}%
\usepackage{booktabs}%
\usepackage{algorithm}%
\usepackage{algorithmicx}%
\usepackage{algpseudocode}%
\usepackage{listings}%
\usepackage{tabularx}
\usepackage{multirow}
\usepackage{hyperref}

\newcommand{\R}{\mathbb{R}}
\newcommand{\zero}{\mathbb{O}}

\newtheorem{remark}{Remark}%

\begin{document}

\begin{center}
		{\large\bf{Difference of Convex (DC) approach for neural network approximation with uniform loss function}}
	\end{center}
	
	\begin{center}
		{Vinesha Peiris, Nadezda Sukhorukova}
	\end{center}

\abstract{Neural networks (NNs) can be viewed as approximation tools. Traditionally, NNs are relying on gradient and stochastic gradient (SG) methods. There are a number of available computational packages for constructing least squares approximations, while uniform (minimax) approximations are hard due to their nonsmooth nature. It was recently demonstrated that a difference convex (DC) programming approach is an efficient alternative optimiser for NNs. In this paper, we demonstrate that a DC programming approach is also efficient for minimax approximation. In our numerical experiments, we compare a DC-programming approach and ADAMAX, a commonly used method for minimax NN approximations.

}

{\bf Mathematics Subject Classification (2020): }{90C26, 65K05, 41A50, 90C47, 90C90
}\\

{\bf Keywords:}
minimax, Chebyshev approximation, DC-programming, convex analysis, neural networks

\section{Introduction}\label{sec:introduction}
Neural Networks (NNs) are powerful Artificial Intelligence (AI) techniques. The theory of NNs is based on solid mathematical background, established in~\cite{Cybenko,Hornik1991,LeshnoPinkus1993,pinkus_1999}. These works lead to the work of A.~Kolmogorov and V.~Arnold. Their  celebrated Kolmogorov-Arnold Theorem~\cite{Kol57, Arnold57} was an attempt to solve the 13th Hilbert's problem.

Essentially, NNs provide accurate approximations. These approximations are found as solutions of optimisation problems, where the objective functions represent the approximation error. An excellent overview of optimisation techniques for NNs can be found in~\cite{Sun2020OptimDeepLearning}. 

The most popular choice for the loss function is least squares due to their smoothness. The numerical methods are based on gradient-type techniques, for example, Stochastic Gradient (SG). In the case of uniform approximation, the corresponding optimisation problems are nonsmooth and the applicability of SG is questionable.

In~\cite{DCasalternativetoNN} the authors demonstrated that Difference of Convex (DC) algorithms are an efficient alternative to the traditional SG techniques in the case of least squares networks. More details on DC-programming are provided in section~\ref{sec:DC}. 

Despite a significant progress in NN approximation in least squares criterion, its Chebyshev (uniform) counterpart is still behind due to its complexity: the corresponding optimisation problems are nonsmooth and non-convex. Chebyshev approximation is a fascinating mathematical problem itself and moreover, it has several practical applications. For example, in~\cite{PeirisRoshchinaSukh} the authors demonstrated that uniform approximation is a more reliable approach if the amount of data is limited. A comprehensive survey on necessary and sufficient optimality conditions for Chebyshev approximation (with a focus on cone constraints) can be found in~\cite{DolgopolikUnified}. See also references therein.

The contribution of this paper is as follows. First of all, we formulate uniform-based NN as DC-programming problems. Then, we implement the numerical experiments with DC-algorithms. Finally, we compare the results with the performance of ADAMAX~\cite{diederik2014adam}, a popular technique for uniform networks.

\section{Motivation}\label{sec:motivation}
 
DC-programming is a powerful tool in the area of modern optimisation, namely, global non-convex optimisation. The origins of DC-programming appear in~\cite{TAO1986249}. Then there have been decades of fruitful development of the theory and algorithms (DCA, stands for DC-programming Algorithm). DCA and its variations have been successfully applied to real-world problems~\cite{An2005,Strekalovsky2017} just to name a few. A detailed modern review of DC-programming and the corresponding algorithms can be found in~\cite{LeThi2018,OpenIssuesDCA}.

One of the applications, where DC-programming and DCA have been successfully implemented is deep learning and neural networks. In particular, in~\cite{Awasthi2024} the authors use DCA for the computation of the adversarial attack for adversarial
robustness. Also, in~\cite{DCasalternativetoNN} the authors demonstrated that DCA can be used as an alternative approach to the neural networks. These papers demonstrated the efficiency of DC-programming in this application. 

In the current paper, we extend the existing results in the following directions.
\begin{enumerate}
    \item Instead of least squares approximation, we use uniform ($\max$-norm) based loss function. This approach may be more efficient in the case when the available data are limited~\cite{PeirisRoshchinaSukh}.
    \item We use ReLU as the activation function. This function is commonly used in neural networks, but its nonsmooth nature limits the application of gradient-based optimisation techniques, in particular, when the loss function is nonsmooth, which is the case for uniform approximation.
\end{enumerate}

It is mentioned in~\cite{OpenIssuesDCA} that in the case of ReLU and Leaky ReLU activation functions, the implementation of DCA can be reduced to solving linear programming problems repeatedly. This is a great advantage, since modern optimisation packages are able to solve large-scaled linear programming problems efficiently. In this paper, we explicitly formulate these linear programming problems and solve them within DCA. 

In our numerical experiments, we compare the results obtained by DCA with ADAMAX~\cite{diederik2014adam}, a specially developed algorithm for uniform-based loss functions (minimax approximation).

In the next section, we provide all the necessary mathematical constructions: DC decomposition (that is, the presentation as a difference of convex functions) of the networks and the implementation of DCA.

\section{DC decomposition and DCA implementation}\label{sec:DC}
\subsection{Neural network approximations}
Deep learning and neural networks are popular tools for function and data approximation. In this paper, we concentrate on the purely approximation nature of the problem.  

It was demonstrated in~\cite{LeshnoPinkus1993,pinkus_1999} that even neural networks with a single hidden layer are efficient approximation techniques. In most practical problems, the goal is to reconstruct (approximate) a function, whose values are known in a finite number of points (dataset $Q$).  Essentially, a single hidden layer neural network with ReLU activation function construct functions and data approximations in the following form:
\begin{equation}\label{eq:NNapp}
\sum_{i=1}^{n}\alpha_i\max\{0,{\bf a}_i^T{\bf T}\},
\end{equation}
where ${\bf T}\in Q$ are the discrete points, ${\bf T}\in\R^d$, $\alpha_i$, $i=1,\dots,n$ and the components of vectors, ${\bf a}_i^T\in\R^d$ are the parameters of the neural network. The components of ${\bf a}_i^T\in\R^d$ are also called weight. The number of components $n$ in the sum of~(\ref{eq:NNapp}) is the number of nodes in the hidden layer. In most practical problems, this number is not known. When $n$ increases, the accuracy of the approximations improves, while the computational complexity increases. 

Note that the activation function ReLU~$\sigma(x)=\max\{0,x\}$ is positively homogeneous of degree one. Therefore, if the sign of each $\alpha_i$, $i=1,\dots,n$ is known, the absolute value of $\alpha_i$, $i=1,\dots,n$ can be reconstructed through the components of vectors~${\bf a}_i^T$. Hence, the approximation in~(\ref{eq:NNapp}) can be constructed as follows:
\begin{equation}\label{eq:NNpair}
    \sum_{i=1}^{n_1}\max\{0,{\bf a}_i^T{\bf T}\}-\sum_{i=n_1+1}^{n_1+n_2}\max\{0,{\bf a}_i^T{\bf T}\},
\end{equation}
where $n_1$ and $n_2$ are large enough. Similar decompositions have been proposed in~\cite{TwoApproaches}. If $n_1=n_2$ in~(\ref{eq:NNpair}), one can talk about the number of pairs (one component comes with a plus and another with a minus).

Therefore, the minimax optimisation problem studied in this paper is as follows:
\begin{equation}\label{eq:optprob}
    \text{minimise} \max_{i=1,\dots,N} \left |f({\bf T}_i)- \left (\sum_{j=1}^{n_1}\max\{0,{\bf a}_j^T{\bf T}_i\}-\sum_{j=n_1+1}^{n_1+n_2}\max\{0,{\bf a}_j^T{\bf T}_i\} \right ) \right |.
\end{equation}
In this paper, for simplicity, we assume that $n_1=n_2=n$ and denote the weights from the second group (for $j\geq n_1+1$) by ${\bf b}_j$, $j=1,\dots,n$. Then problem~(\ref{eq:optprob}) is
\begin{equation}
 \label{eq:optprob_n}
    \text{minimise} \max_{i=1,\dots,N} \left |f({\bf T}_i)- \left (\sum_{j=1}^{n}\max\{0,{\bf a}_j^T{\bf T}_i\}-\sum_{j=1}^{n}\max\{0,{\bf b}_j^T{\bf T}_i\} \right ) \right |.   
\end{equation}

In the next subsection, we formulate a DC decomposition to the objective function in~(\ref{eq:optprob_n}).
\subsection{DC decomposition}
Let $p(x)=\max\{f_i(x), i=1,\dots,n\}$, where $f_i$, $i=1,\dots,n$ are DC functions
$$f_i=g_i-h_i,~i=1,\dots,n$$
and $g_i$, $h_i$, $i=1,\dots,n$ are convex functions. A standard approach for DC decomposition of the maximum of DC functions (see, for example,~\cite{OpenIssuesDCA}) is as follows:
\begin{equation}\label{eq:DC}
    f(x)=\max_{i=1,\dots,n}\biggl \{g_i(x)+\sum_{j=1, j\neq i}^{n}h_{j}\biggr \}-\sum_{j=1}^{n}h_j(x)=g(x)-h(x),
\end{equation}
where
$g(x)=\max_{i=1,\dots,n}\{g_i(x)+\sum_{j=1, j\neq i}^{n}h_{j}(x)\}$ 
and $h(x)=\sum_{j=1}^{n}h_j(x)$ are convex functions.

For simplicity, in this paper we assume that in~(\ref{eq:optprob}), $n_1=n_2$ work with the following optimisation problem:
\begin{equation}\label{eq:optpairs}
 \text{minimise} \max_{i=1,\dots,N}
 \left|f({\bf T}_i)-\left(\sum_{j=1}^{n}\max\{0,{\bf a}_j^T{\bf T}_i\}-\sum_{j=1}^{n}\max\{0,{\bf b}_j^T{\bf T}_i\}\right)\right|,    
\end{equation}
subject to ${\bf A}$ and ${\bf B}$, where ${\bf A}$ and ${\bf B}$ are the components of the weight vectors ${\bf a}_j$ and ${\bf b}_j$, $j=1,\dots,n$.

Denote by $\tilde{p}_i({\bf A},{\bf B})$, $i=1,\dots,n$ the expression under and therefore, for each $i=1,\dots,N$, the expression under the absolute value in~(\ref{eq:optpairs}). Then the objective function in~(\ref{eq:optpairs}) is
\begin{equation*}
    p({\bf A},{\bf B})=\max_{i=1,\dots,N}\max\left\{p_i({\bf A},{\bf B}),-p_i({\bf A},{\bf B})\right\}.
\end{equation*}
The function ${p}_i({\bf A},{\bf B})=\max\left\{\tilde{p}_i({\bf A},{\bf B}),-\tilde{p}_i({\bf A},{\bf B})\right\}$ is DC as the maximum of two DC functions. Applying~(\ref{eq:DC}), obtain
\begin{equation*}
    p_i({\bf A},{\bf B})=g_i({\bf A},{\bf B})-h_i({\bf A},{\bf B}),
    \end{equation*}
where     
    $$    g_i({\bf A},{\bf B})= \max \left \{f({\bf T}_i)+2\sum_{j=1}^{n}\max\{0,{\bf b}_j^T{\bf T}_i\},2\sum_{j=1}^{n}\max\{0,{\bf a}_j^T{\bf T}_i\}-f({\bf T}_i)\ \right \}$$
    and
    $$h_i({\bf A},{\bf B})= \sum_{j=1}^{n}\max\{0,{\bf a}_j^T{\bf T}_i\}+\sum_{j=1}^{n}\max\{0,{\bf b}_j^T{\bf T}_i\}.
$$

Finally, the function $$p({\bf A},{\bf B})=\max_{i=1,\dots,N}p_i({\bf A},{\bf B})$$ is also a DC function. A possible DC decomposition is as follows:
\begin{equation}\label{eq:mainDC}
    p({\bf A},{\bf B})=g({\bf A},{\bf B})-h({\bf A},{\bf B}),
\end{equation}
where
\begin{equation}\label{eq:DCg}
g({\bf A},{\bf B})=\max_{i=1,\dots,N}\left\{g_i({\bf A},{\bf B})+\sum_{j=1,j\neq i}^{N}h_i({\bf A},{\bf B})\right\}\end{equation}
and 
\begin{equation}\label{eq:DCh}h({\bf A},{\bf B})=\sum_{j=1}^{N}h_i({\bf A},{\bf B}).
\end{equation}

To summarise, we have a DC decomposition of the objective function. In the next subsection, we provide the description of DCA, developed in~\cite{TAO1986249}. 

\subsection{DCA}
In our experiments, we use Algorithm~1, also known as Generic DCA or just DCA. There are several modifications and improvements of this method. For more information, see \cite{LeThi2018,OpenIssuesDCA} and references therein. 
\begin{remark}
It is important to note that DCA is the own name of Algorithm~1 and therefore the name DCA should not be used for other algorithms for DC programs.
\end{remark}

\begin{algorithm}
\caption{Algorithm 1: Generic DCA}
\begin{algorithmic}[1] 
\item[] {\bf Initialisation:} DC decomposition: $f(x)=g(x)-h(x)$, initial point $x_0$, $k=0$
\item[]{\bf Repeat}
\item[{\bf Step 1}]     Compute $y_k\in \partial h(x_k)$ where $\partial h(x_k)$ is the subdifferential of $h$ at $x_k$.

\item[{\bf Step 2}]     Compute $x_{k+1}=\arg\min_{x\in\R^n}\left\{g(x)-\langle y_k,x\rangle\right\}$
\item[{\bf Step 3}]    Set $k\leftarrow k+1$.
\item[]{\bf Until}: Stopping criterion
\end{algorithmic}
\end{algorithm}
In Algorithm~1, one replaces the minimisation of a nonconvex function with a sequence of convex optimisation problems. In the next subsection, we demonstrate that the implementation of DCA straightforward:
\begin{itemize}
    \item the subdifferential $\partial h$ is a polytope or even a singleton, whose computation is clear (DCA: Step~1);
    \item the optimisation problem (DCA: Step~2) can be formulated as a linear programming problem.
\end{itemize}

\subsection{Subgradient of \ensuremath{h}: DCA Step 1}

DCA requires the computation of a subgradient of function~$h$ in the decomposition~(\ref{eq:DC}). Recall that ${\bf a}_i, {\bf b}_i\in\R^d$, $i=1,\dots,n$ are the weights and also the decision variables. In the rest of the paper, we denote the decision variables (weights) by 
\begin{equation}\label{eq:weight_vector}
{\bf w}=\begin{bmatrix}
    {\bf A}\\
    {\bf B}
\end{bmatrix},\quad {\bf W}\in\R^{2nd},
\end{equation}
where ${\bf A}\in \R^{nd}$ and ${\bf B}\in\R^{nd}$ are the column vectors that combine the weights ${\bf a}_i$ and ${\bf b}_i$ respectively, $i=1,\dots,n$.

Function~$h$ is the sum of convex piecewise-linear functions. The subdifferentials of each individual function~$h_i$ is a singleton or a polytope. Therefore, the subdifferential of $h$ is also a singleton (if the subfifferentials of all the functions $h_i$ are singletons) or a polytope (general case). 

 For each $j=1,\dots,n$, compute the row vectors~${\bf y}_j^+$ and ${\bf y}_j^-$, $j=1,\dots,n$, such that 
    $${\bf y}_j^+=\zero_n+\sum_{i=1,\\ {\bf a}^T_j{\bf T}_i>0}^{N}{\bf T}_i \quad \text{and}\quad {\bf y}_j^-=\zero_n+\sum_{i=1,\\ {\bf b}^T_j{\bf T}_i>0}^{N}{\bf T}_i.$$
    Then vector
    $${\bf y}=
    \begin{bmatrix}
    {\bf y}_1^+\\{\bf y}_2^+\\
    \vdots \\{\bf y}_n^+\\{\bf y}_1^-\\{\bf y}_2^-\\
    \vdots\\
    {\bf y}_n^- 
    \end{bmatrix}
    $$
is a subgradient of $h$ at ${\bf w}$, where the vector ${\bf w}$ is the collection of weights as in~(\ref{eq:weight_vector}).
\subsection{Uniform approximation}\label{subsec:uniform}

\subsubsection{Minimisation of the linearised surrogate: DCA Step 2, ReLU activation}

The optimisation problem in Step~3 of Algorithm~1 is convex. We can demonstrate that in our specific problem, this problem can be reformulated as a linear programming problem. Indeed, the problem in Step~2 of Algorithm~1 can be formulated as follows:
\begin{equation}\label{eq:z1}
    \text{minimise} \quad z
\end{equation}
subject to
\begin{equation}\label{eq:con1}
 f({\bf T}_i)+2\sum_{j=1}^{n}\max\{0,{\bf b}_j^T{\bf T}_i\}-{\bf y}^T{\bf w}   \leq z,\quad i=1,\dots,N,\\
 \end{equation}
 \begin{equation}\label{eq:con2}
 2\sum_{j=1}^{n}\max\{0,{\bf a}_j^T{\bf T}_i\}-f({\bf T}_i)-{\bf y}^T{\bf w} \leq z,\quad i=1,\dots,N.    
\end{equation}
This problem can be reformulated as a linear programming problem. We introduce new variables:
\begin{equation}
    z_{ij}^+=\max\{0,{\bf a}_j^T{\bf T}_i\},\quad z_{ij}^-=\max\{0,{\bf b}_j^T{\bf T}_i\}, 
\end{equation}
then the optimisation problem~(\ref{eq:z1})-(\ref{eq:con2}) is as follows:
\begin{equation}\label{eq:z2}
    \text{minimise} \quad z
\end{equation}
subject to
\begin{equation}\label{eq:con12}
 f({\bf T}_i)+2\sum_{j=1}^{n}z_{ij}^--{\bf y}^T{\bf w}   \leq z,\quad i=1,\dots,N,\\
 \end{equation}
 \begin{equation}\label{eq:con22}
 2\sum_{j=1}^{n}z_{ij}^+-f({\bf T}_i)-{\bf y}^T{\bf w} \leq z,\quad i=1,\dots,N,    
\end{equation}
\begin{equation}\label{eq:con32}
    z_{ij}^+\geq 0,\quad z_{ij}^+\geq {\bf a}_j^T{\bf T}_i,\quad i=1,\dots,N,\quad j=1,\dots,n,
    \end{equation}
    \begin{equation}\label{eq:con42}
    z_{ij}^-\geq 0,\quad z_{ij}^-\geq {\bf b}_j^T{\bf T}_i,\quad i=1,\dots,N,\quad j=1,\dots,n.
    \end{equation}
    Problem~(\ref{eq:z2})-(\ref{eq:con42}) is a linear programming problem and therefore can be solved efficiently. The total number of decision variables is $2dn+1+2Nn$.
\subsubsection{Minimisation of the linearised surrogate: DCA Step 2, Leaky ReLU activation} 

In the case of Leaky ReLU activation function, the procedure is similar. In this case, the activation function is $\sigma(x)=\max\{\alpha x,x\}$, where $\alpha$ is a small positive number ($\alpha=0.01$ is a common choice). 

Similar to ReLU, the subgradient $\partial h$ is a polytope or a singleton (DCA, Step~1) and can be computed as follows. For each $j=1,\dots,n$, compute the row vectors~${\bf y}_j^+$ and ${\bf y}_j^-$, $j=1,\dots,n$, such that 
    $${\bf y}_j^+=\sum_{i=1,\\ {\bf a}^T_j{\bf T}_i>0}^{N}{\bf T}_i+\alpha\sum_{i=1,\\ {\bf a}^T_j{\bf T}_i\leq 0}^{N}{\bf T}_i  \quad \text{and}\quad {\bf y}_j^-=\sum_{i=1,\\ {\bf b}^T_j{\bf T}_i>0}^{N}{\bf T}_i+ \alpha\sum_{i=1,\\ {\bf b}^T_j{\bf T}_i\leq 0}^{N}{\bf T}_i.$$
    Then vector
    $${\bf y}=
    \begin{bmatrix}
    {\bf y}_1^+\\{\bf y}_2^+\\
    \vdots \\{\bf y}_n^+\\{\bf y}_1^-\\{\bf y}_2^-\\
    \vdots\\
    {\bf y}_n^- 
    \end{bmatrix}
    $$
is a subgradient of $h$ at ${\bf w}$, where the vector ${\bf w}$ is the collection of weights.

The optimisations problem (DCA: Step~2) is an LPP: 
\begin{equation*}
    \text{minimise} \quad z
\end{equation*}
subject to
\begin{equation*}
 f({\bf T}_i)+2\sum_{j=1}^{n}z_{ij}^--{\bf y}^T{\bf w}   \leq z,\quad i=1,\dots,N,\\
 \end{equation*}
 \begin{equation*}
 2\sum_{j=1}^{n}z_{ij}^+-f({\bf T}_i)-{\bf y}^T{\bf w} \leq z,\quad i=1,\dots,N,    
\end{equation*}
\begin{equation*}
    z_{ij}^+\geq \alpha{\bf a}_j^T{\bf T}_i,\quad z_{ij}^+\geq {\bf a}_j^T{\bf T}_i,\quad i=1,\dots,N,\quad j=1,\dots,n,
    \end{equation*}
    \begin{equation*}
    z_{ij}^-\geq \alpha{\bf b}_j^T{\bf T}_i,\quad z_{ij}^-\geq {\bf b}_j^T{\bf T}_i,\quad i=1,\dots,N,\quad j=1,\dots,n.
    \end{equation*}
    
\subsection{Manhattan norm~\texorpdfstring{$L_1$}{L1}}

The subgradient is computed in the same way as it is in subsection~\ref{subsec:uniform} for ReLU and Leaky ReLU activation respectively. The optimisation problems (Step~2 of Algoriyhm~1) can be formulated as LPPs. Later in this section, we demonstrate that the corresponding LPPs in the case of $L_1$ have the same number of constraints as in the case of the uniform approximation, while the number of variables is larger in the case of $L_1$.
\subsubsection{Minimisation of the linearised surrogate: DCA Step 2, ReLU activation}
\begin{equation*}
    \text{minimise} \quad \sum_{i=1}^{N}z_i
\end{equation*}
subject to
\begin{equation*}
 f({\bf T}_i)+2\sum_{j=1}^{n}z_{ij}^--{\bf y}^T{\bf w}   \leq z,\quad i=1,\dots,N,\\
 \end{equation*}
 \begin{equation*}
 2\sum_{j=1}^{n}z_{ij}^+-f({\bf T}_i)-{\bf y}^T{\bf w} \leq z,\quad i=1,\dots,N,    
\end{equation*}
\begin{equation*}
    z_{ij}^+\geq \alpha{\bf a}_j^T{\bf T}_i,\quad z_{ij}^+\geq {\bf a}_j^T{\bf T}_i,\quad i=1,\dots,N,\quad j=1,\dots,n,
    \end{equation*}
    \begin{equation*}
    z_{ij}^-\geq \alpha{\bf b}_j^T{\bf T}_i,\quad z_{ij}^-\geq {\bf b}_j^T{\bf T}_i,\quad i=1,\dots,N,\quad j=1,\dots,n.
    \end{equation*}

\subsubsection{Minimisation of the linearised surrogate: DCA Step 2, Leaky ReLU activation}

The LPP from Step~2 is as follows:
\begin{equation*}
    \text{minimise} \quad \sum_{i=1}^{N}z_i
\end{equation*}
subject to
\begin{equation*}
 f({\bf T}_i)+2\sum_{j=1}^{n}z_{ij}^--{\bf y}^T{\bf w}   \leq z_i,\quad i=1,\dots,N,\\
 \end{equation*}
 \begin{equation*}
 2\sum_{j=1}^{n}z_{ij}^+-f({\bf T}_i)-{\bf y}^T{\bf w} \leq z_i,\quad i=1,\dots,N,    
\end{equation*}
\begin{equation*}
    z_{ij}^+\geq \alpha{\bf a}_j^T{\bf T}_i,\quad z_{ij}^+\geq {\bf a}_j^T{\bf T}_i,\quad i=1,\dots,N,\quad j=1,\dots,n,
    \end{equation*}
    \begin{equation*}
    z_{ij}^-\geq \alpha{\bf b}_j^T{\bf T}_i,\quad z_{ij}^-\geq {\bf b}_j^T{\bf T}_i,\quad i=1,\dots,N,\quad j=1,\dots,n.
    \end{equation*}

Similar to ReLU activation, instead of one variable $z$, we now have $N$ variables $z_i$, $i=1,\dots N$.

\section{Experiments}\label{sec:experimants}
\subsection{Experimental setup}
In this section, we approximate four functions, defined on a finite grid. The first two functions correspond to a dataset with two classes. Essentially, these two functions can only take two values: zero or one. The functions themselves represent the training and test set. The last two functions are uniform discretisations (same discretisation step along each direction) of two continuous functions.
\begin{enumerate}
    \item Two Lead ECG dataset (training set) from~\cite{UCRArchive2018}. This dataset contains 370~signal segments, each signal contains 83~recordings. The output (function value) is binary $\{0,1\}$. The goal is to separate between these two classes.
    \item Two Lead ECG dataset (test set) from~\cite{UCRArchive2018}. This dataset contains 1000~signal segments, each signal contains 83~recordings. The output (function value) is binary $\{0,1\}$. The goal is to separate between these two classes.
    \item Discretisation of $\Phi_1(x,y)=\sqrt{|x-0.5|+3|y|}$ (function with a deep minimum) on the hypercube, where $(x,y) \in [-1,1]\times[-1,1]$. 
    The total number of discretisation points is~2500, the discretisation step is~0.05 along each direction. 
 \item   $\Phi_2(x,y)=\sin(5x-0.5)-\sqrt{|\cos(7y)|}$ (function with several shallow local minima) on the hypercube, where $(x,y) \in [-1,1]\times[-1,1]$. 
    The total number of discretisation points is~2500, the discretisation step is~0.05 along each direction.     
\end{enumerate}
\subsection{Results}

In this section, we present the results of numerical experiments. We compare the results obtained by neural networks and DCA. 

In the case of neural networks, we use~ADAMAX for uniform approximation, while for the Manhattan norm we use ADAM. In all the experiments with neural networks, the computational time was just a fraction of a second. The number of internal nodes (single hidden layer) is either~2 or~4. 

In the case of DCA, the computational time is significantly higher: up to 30~minutes. ``F'' corresponds to the case when the program was running for more than 30~minutes. ``Zero'' corresponds to the objective function being zero (global minimum): interpolation. The number of pairs is one or two, which roughly corresponds to 2 and 4~nodes.

\begin{table}[]
    \centering
    \begin{tabular}{|c|c|c|c|c|c|}
    \hline
Function  &Activation& Loss  & Number of  & Objective &Objective\\
          & function & function & internal& function &  function\\
         &          &          & nodes/pairs  &  value (NN) & value (DCA)\\
\hline
Two Lead ECG & \multirow{2}{*}{ReLU}  & \multirow{2}{*}{Uniform} 
             & 2 (1 pair) & 0.502378  & Zero\\
training     &            &           & 4 (2 pairs) & 0.517571 & zero\\
    \hline
Two Lead ECG & \multirow{2}{*}{Leaky ReLU} & \multirow{2}{*}{Uniform}  
             & 2 (1 pair) & 0.485107  & Zero\\
training     &            &           & 4 (2 pairs) & 0.500470 & zero\\
    \hline
Two Lead ECG & \multirow{2}{*}{ReLU}  & \multirow{2}{*}{Manhattan}
             & 2 (1 pair) & 51.82849  & Zero\\
training     &            &           & 4 (2 pairs) & 84.64083 & zero\\
    \hline
Two Lead ECG & \multirow{2}{*}{Leaky ReLU} & \multirow{2}{*}{Manhattan} 
             & 2 (1 pair) & 32.75499  & Zero\\
training     &            &           & 4 (2 pairs) & 84.97383 & zero\\
\hline
\hline
Two Lead ECG & \multirow{2}{*}{ReLU}  & \multirow{2}{*}{Uniform} 
             & 2 (1 pair) & 0.570136  & F\\
 test        &            &           & 4 (2 pairs) & 0.538044 & F\\
     \hline
Two Lead ECG & \multirow{2}{*}{Leaky ReLU} &\multirow{2}{*}{Uniform} 
             & 2 (1 pair) & 0.503036  & 0.3426\\
test         &            &           & 4 (2 pairs) & 0.547249 & F\\
    \hline
Two Lead ECG & \multirow{2}{*}{ReLU}  & \multirow{2}{*}{Manhattan} 
             & 2 (1 pair) & 98.170    & F\\
 test        &            &           & 4 (2 pairs) & 122.974  & F\\
     \hline
Two Lead ECG & \multirow{2}{*}{Leaky ReLU} & \multirow{2}{*}{Manhattan} 
             & 2 (1 pair) & 135.594   & 86.4473\\
test         &            &           & 4 (2 pairs) & 124.185  & F\\
\hline
\hline
$\Phi_1$    & \multirow{2}{*}{ReLU}   & \multirow{2}{*}{Uniform} 
            & 2 (1 pair)  & 0.923042  & 0.8835\\
            &             &           & 4 (2 pairs) & 0.876730 & 0.8813\\
     \hline
$\Phi_1$    & \multirow{2}{*}{Leaky ReLU} & \multirow{2}{*}{Uniform} 
            & 2 (1 pair)  & 0.950048  & 0.8835\\
            &             &           & 4 (2 pairs) & 0.883299 & 0.8830\\
\hline
$\Phi_1$    & \multirow{2}{*}{ReLU}   & \multirow{2}{*}{Manhattan}
            & 2 (1 pair)  & 717.0725  & 691.8513\\
            &             &           & 4 (2 pairs) & 325.435  & 137.8907\\
     \hline
$\Phi_1$    & \multirow{2}{*}{Leaky ReLU} & \multirow{2}{*}{Manhattan} 
            & 2 (1 pair)  & 719.205   & 691.8513\\
            &             &           & 4 (2 pairs)  & 331.3275 & 141.2259\\
\hline
\hline
$\Phi_2$    & \multirow{2}{*}{ReLU}   & \multirow{2}{*}{Uniform} 
            & 2 (1 pair)  & 1.630845  & 1.9937\\
            &             &           & 4 (2 pairs)  & 1.424211  & 1.9937\\
     \hline
$\Phi_2$    & \multirow{2}{*}{Leaky ReLU} & \multirow{2}{*}{Uniform} 
            & 2 (1 pair)  & 1.630653  & 1.9937\\
            &             &           & 4 (2 pairs)  & 1.409082  & 1.9937\\
\hline
$\Phi_2$    & \multirow{2}{*}{ReLU}   & \multirow{2}{*}{Manhattan} 
            & 2 (1 pair)  & 1640.4725 & 1361.5\\
            &             &           & 4 (2 pairs)  & 1581.4775 & 1361.5\\
     \hline
$\Phi_2$    & \multirow{2}{*}{Leaky ReLU} & \multirow{2}{*}{Manhattan} 
            & 2 (1 pair)  & 1637.65   & 1361.5\\
            &             &           & 4 (2 pairs)  & 1553.265  & 1361.5\\
\hline
    \end{tabular}
    \caption{Numerical experiments: neural networks (NN) and DCA}
    \label{table:NN_DCA}
\end{table}

The results are in Table~\ref{table:NN_DCA}. The last two columns correspond to the final value of the objective function: column NN corresponds to neural networks and DCA corresponds to DC-optimisation. From this table, one can see that DCA performs significantly better that neural networks on Two Lead ECG training (smaller dataset), where DCA reached a global solution (interpolation). In the case of the larger dataset (Two Lead ECG test), DCA failed to produce results within 30~minutes, but in the cases when the solution was achieved (Leaky ReLU, one pair, uniform and Manhattan distance) the objective function value is reached by DCA is significantly better. In the case of $\Phi_1$, DCA works better than neural networks (except ReLU, uniform, 2 pairs, where the results are close to each other, but neural networks are slightly better). In the case of $\Phi_2$, DCA works better for Manhattan  approximation, while neural networks are more efficient for uniform approximation.

\section{Conclusions and future research}\label{sec:conclusions}

In this paper, we demonstrated that the optimisation problems in neural networks approximation (single hidden layer) for uniform and Manhattan approximation can be formulated and solved using DC-optimisation. The numerical experiments demonstrate that the computational time is significantly higher for DCA (all experiments), while the objective function values are better in the case of DCA (on smaller size problems, DCA reached a globally optimal solution).

In the future, we will be working on the improvement of the computational time  for DCA. In particular, instead of solving a large linear programming problem, it may be reverted to a linear search. One of the options is to use a line search from~\cite{BacktrackingGD}.  




\section*{Declarations}

\begin{itemize}
\item Funding: Not applicable.
\item Conflict of interest/Competing interests:  Not applicable.
\item Ethics approval and consent to participate: Not applicable.
\item Consent for publication: All authors give their consent for publication
\item Data availability: Data are publicly available. References are provided. 
\item Materials availability: Not applicable.
\item Code availability: Yes, references are provided. 
\item Author contribution: Equal contribution.
\end{itemize}

\noindent










\end{document}